\documentclass[12pt, a4paper]{article}

\usepackage{amsfonts}

\newcommand{\dis}{\displaystyle}

\title{Frobenius $n$-homomorphisms, transfers and branched coverings}
\author{V.~M.~Buchstaber and E.~G.~Rees}
\begin{document}

\maketitle

\noindent {\bf Abstract} The main purpose is to characterise
continuous maps that are $n$-branched coverings in terms of induced
maps on the rings of functions.  The special properties of Frobenius
$n$-homomorphisms between two function spaces that correspond to
$n$-branched coverings are determined completely.  Several
equivalent definitions of a Frobenius $n$-homomorphism are compared
and some of their properties are proved. An axiomatic treatment of
$n$-transfers is given and properties of $n$-branched coverings are
studied and compared with those of regular coverings.

\vspace{0.7cm} \noindent {\bf Introduction} \medskip

In previous papers  \cite{[BR1]}, \cite{[BR2]}, and  \cite{[BR3]} we
have studied the relationship between  Frobenius $n$-homomorphisms
and symmetric products, in particular we showed that  Frobenius
$n$-homomorphisms $C(X) \to C(Y)$ correspond precisely to continuous
maps $f: Y \to {\rm Sym}^n(X)$, thus generalising the classical
theorem of Gelfand-Kolmogorov. In this paper we characterise
$n$-branched coverings using transfer maps which we treat as
Frobenius $n$-homomorphisms of a special type. In \S1 we review and
prove the equivalence of a number of definitions of Frobenius
$n$-homomorphisms that were not covered in our previous papers. \S 2
is devoted to showing that the composition of a Frobenius
$n$-homomorphism with a Frobenius $m$-homomorphism is a Frobenius
$nm$-homomorphism, the proof depends heavily on a previous result
from \cite{[BR1]}. A version of this result appeared in a
preliminary preprint of the present paper and D.~V.~Gugnin had given
an improvement of that result. \S 3 reviews the definition of an
$n$-branched covering due to L.~Smith \cite{[Sm]} and A.~Dold
\cite{[D]} and we prove some new properties and give examples; in
particular we make a detailed analysis of $n$-branched coverings
over an interval. \S 4 introduces the notion of an $n$-transfer for
a ring homomorphism as a special kind of Frobenius $n$-homomorphism.
It is shown that the kernel introduced by Gugnin \cite{[G]} for a
Frobenius $n$-homomorphism is trivial for an $n$-transfer. \S 5
contains the main theorem which characterises $n$-branched coverings
in terms of $n$-transfers between the relevant function spaces.

\vspace{1.3cm}

\noindent \S 1 {\bf Frobenius $n-$homomorphisms}

\bigskip
In \cite{[BR1]}, \cite{[BR2]} we introduced the concept of a
Frobenius $n$-homomorphism and studied some of its properties (these
papers also contain references to related work by other authors).
We recall the basic definition.

Consider a linear map $ f : A \to B$ between two commutative,
associative ${\mathbb C}$ algebras. The maps $\Phi _{n}(f):A^{\otimes
n}\rightarrow B$ are  defined as follows:

Each permutation $\sigma \in \Sigma_{n}$ the symmetric group on $n$
letters, can be decomposed into a product of disjoint cycles of
total length $n,$ say $ \sigma =\gamma _{1}\gamma _{2}\dots \gamma
_{r}.$ If\ $\gamma =(i_{1}\ldots i_{m})$ is a cycle, let $f_{\gamma
}(a_{1},a_2,\ldots ,a_{n})=f(a_{i_{1}}a_{i_{2}}\ldots a_{i_{m}})$
then {\footnotesize $$ \Phi _{n}(f)(a_{1},a_2,\ldots ,a_{n})=
\sum\limits_{\sigma \in \Sigma_{n}} \varepsilon _{\sigma}f_{\gamma
_{1}}(a_{1},a_2,\ldots,a_{n})f_{\gamma _{2}}(a_{1},a_2,\ldots
,a_{n})\ldots f_{\gamma _{r}}(a_{1},a_2,\ldots ,a_{n})$$} where $
\varepsilon _{\sigma}$ is the sign of the permutation $\sigma .$

From this definition it is clear that $\Phi _{n}(f):A^{\otimes
n}\rightarrow B$ is $n$-linear and symmetric, one can use
polarisation to ease the verification of some of its properties; in
other words to prove identities it is enough to consider the values
of $\Phi _{n}(f)(a,a,\ldots,a)$ for all $a \in A$ (we will sometimes
abbreviate this to $\Phi _{n}(f)(a)$).

There is also an inductive definition (first used by Frobenius in
the case of group algebras of finite groups) for the $\Phi _{n}(f)$
starting with $\Phi _{1}(f)=f$ and, for $n \geq 1,$ {\footnotesize
$$\Phi_{n+1}(f)(a_0,a_1,\ldots,a_n) = f(a_0)\Phi
_{n}(f)(a_1,a_2,\ldots,a_n ) - \sum_{r=1}^n \Phi
_{n}(f)(a_1,a_2,\ldots,a_0a_r, \ldots, a_n )$$} or equivalently,
because of the polarisation identities for symmetric multilinear
maps, {\footnotesize
$$\Phi_{n+1}(f)(a,a,\ldots,a)= f(a)\Phi_{n}(f)(a,a,\ldots,a)
 -n \Phi_{n}(f)(a^2,a,\ldots,a).$$}
We will also find it useful
to have another equivalent  definition:

\noindent
$\Phi_n(f)(a)$ is the determinant of the
matrix
$$
   \left(
\begin{array}{cccccc}
f(a) & 1 & 0 & 0& \dots &0 \\
f(a^2) & f(a) & 2  & 0& \dots &0 \\
\vdots &\vdots & \vdots& \ddots& & \vdots \\
\vdots&&&&&\vdots\\
f(a^{n-1}) &f(a^{n-2}) & f(a^{n-3}) & \dots & f(a)
& {n-1 } \\
f(a^n) &f(a^{n-1}) & f(a^{n-2})  &  \dots &f(a^2)
&f(a)

\end{array} \right)
            $$

\bigskip

\noindent {\bf Definition 1.1} \hspace{0.3cm} A linear map $f : A
\to B$ is called a {\bf Frobenius

\noindent
$n$-homomorphism} if it satisfies $f(1)=n$ and $\Phi _{n+1}(f)\equiv 0.$

\medskip We note that, by polarisation, the condition $\Phi _{n+1}(f)\equiv 0$
is equivalent to $\Phi _{n+1}(f)(a) = 0$ for all $a \in A.$

\bigskip
\noindent {\bf Definition 1.2} \hspace{0.3cm} An algebra $A$ is {\bf
connected} if, for each $k \geq 1$, the equation

\noindent
$x(x-1)\ldots (x-k)=0$ has only the obvious $k+1$ solutions.

\smallskip
This terminology is suggested by the fact that an algebra of functions
on a space $X$ has this property if and only if $X$ is connected.

\smallskip
The following result is easily proved by induction

\bigskip
\noindent {\bf Lemma 1.3} \hspace{0.3cm} {\it For any linear map $f$
one has $\Phi_{n+1}(f)(a,1,\ldots,1)= f(a)(f(1)-1)(f(1)-2) \ldots
(f(1)-n).$}

\smallskip
By taking $a=1$ one obtains

\bigskip
\noindent {\bf Corollary  1.4} \hspace{0.3cm} {\it  Let $B$ be
connected and  $f : A \to B$ be such that $\Phi_{n+1}(f)\equiv 0$
then $f(1) \in \{0,1,2,\ldots,n\}.$}

\medskip
These results are closely related to those of Corollary 2.5 of
\cite{[BR1]}.

\medskip
The condition $f(1)=n$ plays a crucial r\^ole (Proposition 2.7 of
\cite{[BR1]}) in the definition of a Frobenius $n$-homomorphism, as
the following indicates.

\bigskip
\noindent {\bf Proposition 1.5} \hspace{0.3cm} {\it If
$\Phi_{n+1}(f)\equiv 0$ and $f(1)=k$ then $\Phi_{k+1}(f)\equiv 0.$}

\medskip
We will use the following result about symmetric polynomials (it is
related to some of those in \cite{[BR3]}, \cite{[BR4]}); it helps to
explain the determinant expression used above.

\bigskip
\noindent {\bf Proposition 1.6} \hspace{0.3cm} { \it Let $s_k =
\beta_1^k + \beta_2^k + \ldots + \beta_n^k$, then the indeterminates
$\beta_r \; (1 \leq r \leq n)$ are the roots of the polynomial
$d(t)$ given as the determinant of the matrix

$$
   \left(
\begin{array}{ccccccc}
s_1 & 1 & 0 & 0& \dots &\dots&0 \\
s_2 & s_1 & 2  & 0& \dots  &\dots&0 \\
\vdots &\vdots & \vdots& \ddots&& & \vdots \\
\vdots &&&&&&\vdots\\
s_{n-1} &s_{n-2} & s_{n-3} & \dots & s_1
& {n-1 }&0 \\
s_n &s_{n-1} & s_{n-2}  &  \dots &s_2
&s_1&n\\
t^n&t^{n-1}& t^{n-2}& \dots&t^2&t&1

\end{array} \right)
            $$

Moreover, if $f : A \to B$ is linear and $f(a^k)=s_k$ then the
determinant $d(t)$  is, up to a non-zero constant multiple, equal to
$$t^n-\Phi_1(f)(a)t^{n-1}+\frac{1}{2}\Phi_2(f)(a)t^{n-2}- \ldots +\frac{(-1)^n}{n}\Phi_n(f)(a).$$}
\noindent
{\bf Proof}

As usual, let $e_r$ denote the elementary symmetric polynomial of
degree $r$ in $\beta_r \; (1 \leq r \leq n)$.
Denote the columns of the above matrix by
 $({\bf c}_1,{\bf c}_2,\ldots, {\bf c}_{n+1})$ and
replace the first column by the column vector $${\bf c}_1 - e_1{\bf
c}_2 + e_2{\bf c}_3 - \ldots +(-1)^ne_n{\bf c}_{n+1}.$$ The
determinant is unchanged.
 Using the standard Newton formulae (see \cite{[Mac]}, page
20),
$$s_r - s_{r-1}e_1 + \ldots +(-1)^{r-1}s_1e_{r-1} +(-1)^r re_r =0$$
\noindent we see that the first column of the new matrix has zero
entries except for the last entry which equals $$ p(t) = t^n
-e_1t^{n-1} +e_2t^{n-2} - \ldots + (-1)^ne_n.$$ So $d(t) = (-1)^n
\,n!\, p(t)$ and hence the roots of $d(t)$ are the same as those of
$ p(t),$ namely $\{\beta_1,\beta_2,\ldots,\beta_n\}.$

Finally, to verify the result about $f : A \to B$ it is enough to
consider the special case $A= \mathbb{C}[a], B=
\mathbb{C}[\beta_1,\beta_2, \ldots ,\beta_n ]$. Then, considering
the appropriate sub-determinant of the above determinant that
defines $e_n,$ we see that $\Phi_n(f)(a)=n!e_n.$

\bigskip
This leads us to yet another characterisation of
Frobenius $n$-homomorphisms (which is closely related to some formulae
in \cite{[BR3]} and \cite{[BR4]}).

\bigskip
\noindent {\bf Proposition 1.7} \hspace{0.3cm} {\it A linear map $f
: A \to B$ is a Frobenius $n$-homomorphism if and only if for each
$a \in A$ there is a polynomial $p_a(t) \in B[t]$ of degree $n$ such
that
$$\sum_{q=0}^{\infty}\frac{f(a^q)}{t^{q+1}} = \frac{{\rm d}}{{\rm d}
t}\log p_a(t).$$}
\noindent {\bf Proof}

Given the Frobenius
$n$-homomorphism  $f : A \to B$ and $a \in A$, consider the polynomial
$\dis{p_a(t)=  \det (M) \in B[t]}$ where   $M$ is the matrix

$$
   \left(
\begin{array}{ccccccc}
f(a) & 1 & 0 & 0& \dots &\dots&0 \\
f(a^2) & f(a) & 2  & 0& \dots  &\dots&0 \\
\vdots &\vdots & \vdots& \ddots&& & \vdots \\
\vdots &&&&&&\vdots\\
f(a^{n-1}) &f(a^{n-2}) & f(a^{n-3}) & \dots & f(a)
& {n-1 }&0 \\
f(a^n) &f(a^{n-1}) & f(a^{n-2})  &  \dots &f(a^2)
&f(a)&n\\
t^n&t^{n-1}& t^{n-2}& \dots&t^2&t&1

\end{array} \right)
            $$

Choose an extension $\overline{B}$ of $B$ such that $p_a(t)$ factors
completely in $\overline{B}[t]$, say
$p_a(t)=n!(t-\beta_1)(t-\beta_2)\ldots (t-\beta_n)$ with $\beta_r
\in \overline{B}.$ Then
$$\frac{{\rm d}}{{\rm d}t}\log p_a(t)=  \frac{{\rm d}}{{\rm
    d}t}\left(\log(t-\beta_1)+ \log(t-\beta_2)+ \ldots + \log(t-\beta_n)\right)$$
But $$ \frac{{\rm d}}{{\rm d}t} \log(t-\beta) = \frac{1}{t-\beta}
=\frac{1}{t}\left(1- \frac{1}{\frac{\beta}{t}}\right)= \frac{1}{t}+
    \frac{\beta}{t^2}+ \frac{\beta^2}{t^3} + \ldots$$
However, by Proposition 1.5, $\beta_1^r +\beta_2^r + \ldots
\beta_n^r
    =f(a^r)$ and the result is proved.

Conversely, if $$\sum_{q=0}^{\infty}\frac{f(a^q)}{t^{q+1}} =
\frac{{\rm d}}{{\rm d} t}\log p_a(t)$$ and
$p_a(t)=n!(t-\beta_1)(t-\beta_2)\ldots (t-\beta_n)$ then
$f(a^r)=\beta_1^r +\beta_2^r + \ldots \beta_n^r$ and, with $r=0$ we
have $f(1)=n$; using this and Proposition 1.5  we get that the
$r^{\rm th}$ elementary symmetric polynomial in $\beta_1, \beta_2,
\ldots, \beta_n$ is $\Phi_r(f)(a)$ and hence that
$\Phi_{n+1}(f)(a)=0.$

\medskip
 As in \cite{[BR1]} we denote the subalgebra of symmetric
tensors in $A^{\otimes n}$ by $\mathcal{S}^nA$.

 We will find the following fact very useful and it is
easy to prove.

\bigskip \noindent {\bf Lemma 1.8} {\it The diagonal map
$\Delta_n :A \to \mathcal{S}^nA$  defined by $$\Delta_n(a) = a
\otimes 1 \otimes \ldots \otimes 1 +1 \otimes a \otimes \ldots
\otimes 1 + \ldots +1 \otimes 1 \otimes \ldots \otimes a$$ is a
Frobenius $n$-homomorphism.}

\bigskip
\noindent {\bf Theorem 1.9} \hspace{0.3cm} {\it  If $f : A \to B$ is
a Frobenius $n$-homomorphism, then the map defined by
$$\frac{\Phi_n(f)}{n!} : \mathcal{S}^nA \to B$$ is a ring
homomorphism. Conversely, if $f : A \to B$ is linear,  $f(1)=n$ and
$\dis{\frac{\Phi_n(f)}{n!}}$ is a ring homomorphism, then $f$ is a
Frobenius $n$-homomorphism.}

\bigskip
\noindent {\bf Proof}

The first statement is Theorem 2.8 of \cite{[BR1]} so we only need
to prove the (easier) converse part.

 Since $\dis{\frac{\Phi_n(f)}{n!}}$ is a ring
homomorphism we have that, for ${\bf a,b} \in \mathcal{S}^nA$
$$\Phi_n(f)({\bf a})\Phi_n(f)({\bf b})= n!\Phi_n(f)({\bf a b}).$$
But by the inductive definition and the symmetry of $\Phi_{n+1}(f)$,
$$\Phi_{n+1}(f)( a^{\otimes n+1})= f(a)\Phi_n(f)( a^{\otimes
n})-n\Phi_n(f)(a^2 \otimes a^{\otimes n-1})$$ and $a^2 \otimes
a^{\otimes n-1}+a \otimes a^2\ldots \otimes a + \ldots + a^{\otimes
n-1}\otimes a^2=a^{\otimes n}\Delta_n(a).$ By Lemma 1.3,
$\Phi_n(f)(\Delta_n(a))=n!f(a).$ So
$$\Phi_{n+1}(f)( a^{\otimes n+1})= f(a)\Phi_n(f)( a^{\otimes n})-
\Phi_n(f)(a^{\otimes n}\Delta_n(a))$$
$$=f(a)\Phi_n(f)( a^{\otimes n})-\Phi_n(f)( a^{\otimes
n})\Phi_n(f)(\Delta_n(a))/n!=\Phi_n(f)( a^{\otimes n})
(f(a)-f(a))=0.$$

\medskip
In an appropriate sense $\Delta_n$ is the universal  Frobenius
$n$-homomorphism on $A$ as the following shows.

\bigskip
\noindent {\bf Corollary 1.10}\hspace{0.3cm} {\it A Frobenius
$n$-homomorphism $f : A \to B$ factors uniquely as  $$\tilde{f}
\Delta_n : A \to \mathcal{S}^nA \to B
$$ where $\tilde{f}$ is a ring homomorphism. Moreover,
$\dis{\tilde{f} = \frac{\Phi_n(f)}{n!}}  .$}

\bigskip
\noindent {\bf Proof}

By Lemma 1.3, $\dis{\frac{\Phi_n(f)}{n!}(\Delta_n(a)) = f(a)}$ which
proves existence.

It remains to prove that $\tilde{f} $ is unique. This is true
because of the fact that, as an algebra, $\mathcal{S}^nA$ is
generated by the elements $\{\Delta_n(a) : a \in A\}$ since a
general element $\sum_{\sigma \in \Sigma_n} a_{\sigma(1)}\otimes
 \ldots \otimes a_{\sigma(n)}$ of $\mathcal{S}^nA$ equals
 $\Phi_n(\Delta_n)(a_1,a_2,\ldots,a_n)$ which is a polynomial in
 elements of the form $\Delta_n(a)$ with $a \in A$.

\medskip
Following an idea due to D.~V.~Gugnin we introduce the appropriate
categorical concept of the kernel of a Frobenius $n$-homomorphism.

\bigskip \noindent
{\bf Definition 1.11} \hspace{0.3cm} For a Frobenius
$n$-homomorphism $f : A \to B$ let $K_f = \{a : f(ax)=0 \;\; {\rm
for \;\;all \;} x \in A\}$, it is clearly a subspace of the kernel
of the linear map $f$.

When $B$ has no nilpotent elements,  D.~V.~Gugnin has shown that
$K_f = \{a : f(a^r)=0  {\rm \;\;for \;\;} 1 \leq r \leq n\}$.

\vspace{0.7cm}
\noindent {\bf \S 2 Compositions}

\medskip
In \cite{[BR1]} we showed that the sum of a Frobenius
$m$-homomorphism and a Frobenius $n$-homomorphism is a Frobenius
$m+n$-homomorphism. The main result of this section is to show that
they also behave appropriately under composition.

\bigskip

\noindent {\bf Theorem 2.1} \hspace{0.3cm} {\it Let $A,B,C$ be
associative, commutative $\mathbb{C}-$algebras.  If $f:A\rightarrow
B $ and $g:B\rightarrow C $ are  Frobenius $n$- and
$m$-homomorphisms, respectively then $gf:A\rightarrow C $ is a
Frobenius $mn$-homomorphism.}

\bigskip\noindent
{\bf Proof}

 We first note that the theorem is trivial if either $n$
or $m$ equals 1 and that we will make constant use of Corollary
1.10.

The map $f$ factors through its universal Frobenius $n$-homomorphism
as $$ \tilde{f} \Delta_n : A \to \mathcal{S}^nA \to B.$$ Similarly,
$g$ factors as $$ \tilde{g} \Delta_m : B \to \mathcal{S}^mB \to C.$$
Since $\tilde{f}$ is a ring homomorphism and $\Delta_m$ is a
Frobenius $m$-homomorphism, the composition
$$\Delta_m \tilde{f} : \mathcal{S}^nA \to B \to \mathcal{S}^mB $$ is
a Frobenius $m$-homomorphism, so it factors through
$$\Delta_m  : \mathcal{S}^nA \to  \mathcal{S}^m\mathcal{S}^nA.$$
By direct observation, the composition
$$\Delta_m \Delta_n : A \to \mathcal{S}^{n}A \to \mathcal{S}^{m}
\mathcal{S}^{n}A$$ is equal to $$i\Delta_{mn} :A \to
\mathcal{S}^{mn}A \to \mathcal{S}^{m}\mathcal{S}^{n}A$$ where $i$ is
the inclusion map. Hence the composition $gf:A \to C$ factors
through $$\Delta_{mn} :A \to \mathcal{S}^{mn}A$$ and so $gf$ is a
Frobenius $mn$-homomorphism.

\vspace{0.7cm}
 \noindent {\bf \S 3 Branched coverings}
\medskip

We consider  branched coverings in the sense studied by Smith
\cite{[Sm]} and by Dold \cite{[D]} and elaborate on their
properties.

 \bigskip
\noindent {\bf Definition 3.1}\hspace{0.3cm} An {\bf $n$-branched
covering} $h:X \to Y$ is a continuous map between two Hausdorff
spaces and a continuous map $t : Y \to {\rm Sym}^n(X)$ such that
\begin{enumerate}
\item[(i)] $x \in th(x)$ for every $x \in X$; and
\item[(ii)] $ {\rm Sym}^n(h)(ty) = ny $ for every $y \in Y.$

\end{enumerate}

\bigskip
\noindent {\bf Example 3.2} \hspace{0.3cm} The map $p: \mathbb{C}
\to \mathbb{C}$ defined by a polynomial of degree $n$ is the classic
example of an $n$-branched covering. The map $t$ is defined by
$t(w)=[z_1,z_2,\ldots,z_n] $ where the $z_r$ are the roots (counted
with multiplicities) of the equation $p(z)=w.$ It is straightforward
to verify the above axioms in this case.

\bigskip
\noindent {\bf Example 3.3} \hspace{0.3cm} If $G$ is a finite group
acting continuously and effectively on a Hausdorff space $X$ and
$h:X \to Y =X/G$ is the map to the space of orbits, then $h$ is an
$n$-branched covering where $n$ is the cardinality of $G$ and the
map $t:Y \to {\rm Sym}^n(X)$ is given by $t(y)$ = the points in the
orbit defined by $y$ (counted with multiplicities). More generally
(see \cite{[D]}, Example 1.4), if $H \subset G$ is a subgroup of
finite index $n$ and $X$ is an effective $G$-space, then the
quotient map $X/H \to X/G$ is an $n$-branched covering.

Every $2$-branched covering arises from an action of the group with
two elements.

\bigskip
\noindent {\bf Proposition 3.4} \hspace{0.3cm} {\it Let $f : X \to
Y, \; s: Y \to {\rm Sym}^n(X)$ and $g : Y \to Z, \; t : Z \to {\rm
Sym}^m(Y)$ be $n$- and $m$-branched coverings, then the composition
$h=gf:X \to Z$ is an $mn$-branched covering with $u : Z \to {\rm
Sym}^{mn}(X)$ being the composition
$$i\;{\rm Sym}^m(s)\;t : Z \to {\rm Sym}^m(Y) \to {\rm Sym}^m({\rm Sym}^n(X)) \to {\rm Sym}^{mn}(X).$$}
The proof is straightforward.

\medskip
One defines induced $n$-branched coverings by taking the obvious
pullback diagram : let $h:X \to Y, \; t: Y\to {\rm Sym}^n(X)$ be an
$n$-branched covering and $\phi : Z \to Y$ a continuous map. The
usual pullback is $\tilde{X} = \{(z,x) | \phi(z)=h(x)\}$, the
projection is $\tilde{h}(z,x)=z$ and $\tilde{t} :Z \to {\rm
Sym}^n(\tilde{X})$ is defined by $
\tilde{t}(z)=[(z,x_1),(z,x_2),\ldots,(z,x_n)]$ where $t
\phi(z)=[x_1,x_2,\ldots, ,x_n].$ The properties of an $n$-branched
covering are straightforward to check.

\medskip
We modify the Ehresmann method \cite{[St]} which  constructs the
principal fibration associated with a locally trivial fibration. As
a result we get a method to `resolve' an $n$-branched covering.

An {\bf epimorphism} from the set $[n]=\{1,2,\ldots,n\}$ to an
$n$-multiset is a map which is onto and the counter image of an
element with multiplicity $m$ has size $m.$

\medskip

 If $h:X \to Y$ is an $n$-branched
covering, let $E$ be the set of all maps $\psi : [n] \to X$ such
that $h\psi(1)=h\psi(2)= \ldots =h\psi(n)$ and $\psi$ is an
epimorphism onto the multiset $th\psi.$ Clearly the symmetric group
$\Sigma_n$ acts on $E$ and the quotient is $Y$; moreover, $E
\times_{\Sigma_n}[n]$ is isomorphic to $X$ and the projection onto
the first factor can then be identified with $h$. The map $\psi$ can
be thought of as a `universal' labelling for the branches of the
covering. This proves the following

\bigskip
\noindent {\bf Theorem 3.5} (\cite{[D]}, Proposition 1.9)
\hspace{0.3cm} {\it Every $n$-branched covering can be described in
the form of a projection $p : E \times_{\Sigma_n}[n] \to
E/{\Sigma_n}$ for some $\Sigma_n$ space $E.$}

\bigskip \noindent
{\bf Example 3.6}

To illustrate all this we consider in some detail the case where the
base space is an interval; this is in contrast to the case of
regular coverings over an interval in which case they are  trivial.
The following discussion describes many $n$-branched coverings of an
interval, including all those with finitely many `branch points'.

 A {\bf set partition }$\pi$ of a finite set $S$ consists of a
family of disjoint non-empty subsets of $S$ whose union is $S$. We
let $n(\pi)$ denote the number of parts of $\pi$.

Let $P_S$ denote the set of all set partitions of $S.$ Two set
partitions $\pi_1=\{A_1,A_2,\ldots,A_p\}, \;
\pi_2=\{B_1,B_2,\ldots,B_q\}$ are called {\bf adjacent} if each
$A_r$ is either a union of some of the $B's$ or is a subset of one
of the $B's$; this relation is symmetric. [For example, if
$S=\{u,v,w\}$ then the two partitions $\{u\},\{v,w\}$ and
$\{u,v\},\{w\}$ are not adjacent but $\{u\},\{v,w\}$ and $\{u,v,w\}$
are adjacent.]
 A map $\phi :[k+1] \to P_{[n]}$ is {\bf compatible} if the
set partitions $\phi(r), \phi(r+1)$ are adjacent for $1 \leq r \leq
k.$

Let  $0=b_0 < b_1 < b_2 < \ldots < b_k < 1=b_{k+1}$ be a dissection
of $I=[0,1].$ Given a compatible $\phi$ we can construct an
$n$-branched covering of $I$ that is branched over  the points $b_r,
\; 1 \leq r \leq k$ as follows :

Over the interval $(b_{r-1},b_r)$ there are $n(\phi(r))$ disjoint
intervals each labelled by a part of the set partition $\phi(r)$. At
the point $b_r$ the two sets of $n(\phi(r))$ and $n(\phi(r+1))$
intervals are joined according to the adjacency between $\phi(r)$
and $\phi(r+1).$

\bigskip
An illustration with $n=5, \; k=2$  and `branched' over two points
$b_1, b_2$ with  $\phi(1) =\{a\},\{b\},\{cde\},$
$\phi(2)=\{a,b\},\{c\},\{d,e\},$ $\phi(3)=\{a,b,c\},\{d\},\{e\}$ is

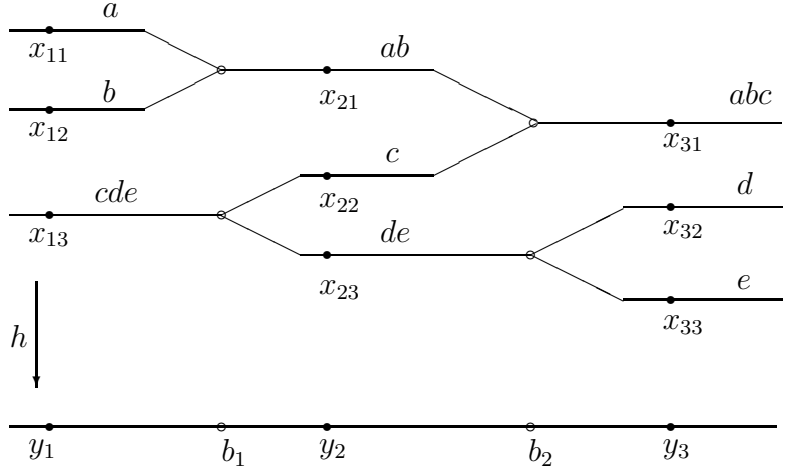
\begin{picture}(200,230)(-30,-120)

\put(0,10){\line(1,0){80}} \put(0,50){\line(1,0){51}}
\put(0,80){\line(1,0){51}}

\put(80,10){\line(2,1){30}} \put(80,10){\line(2,-1){30}}

\put(80,65){\line(-2,1){30}} \put(80,65){\line(-2,-1){30}}

\put(80,65){\line(1,0){80}} \put(110,-5){\line(1,0){87}}
\put(110,25){\line(1,0){50}}


\put(160,25){\line(2,1){39}} \put(160,65){\line(2,-1){39}}

\put(197,-5){\line(2,1){35}} \put(197,-5){\line(2,-1){35}}

\put(200,45){\line(1,0){92}} \put(232,-22){\line(1,0){60}}
\put(232,13){\line(1,0){60}}

\put(0,-70){\line(1,0){290}}

\put(35,85){$a$} \put(35,53){$b$} \put(32,15){$cde$}

\put(140,70){$ab$} \put(142,30){$c$} \put(140,0){$de$}

 \put(10,-15){\vector(0,-1){40}}
\put(0,-40){$h$}

\put(15,10){\circle*{3.0}} \put(15,50){\circle*{3.0}}
\put(15,80){\circle*{3.0}}\put(15,-70){\circle*{3.0}}

\put(7,-80){$y_1$}\put(7,70){$x_{11}$}\put(7,40){$x_{12}$}\put(7,0){$x_{13}$}

\put(120,-5){\circle*{3.0}} \put(120,25){\circle*{3.0}}
\put(120,65){\circle*{3.0}}\put(120,-70){\circle*{3.0}}

\put(117,-80){$y_2$}\put(117,52){$x_{21}$}\put(117,13){$x_{22}$}\put(117,-20){$x_{23}$}

\put(250,45){\circle*{3.0}} \put(250,13){\circle*{3.0}}
\put(250,-22){\circle*{3.0}} \put(250,-70){\circle*{3.0}}

\put(247,-80){$y_3$}\put(247,3){$x_{32}$}\put(247,-33){$x_{33}$}\put(247,36){$x_{31}$}
\put(275,18){$d$}\put(275,-18){$e$}\put(272,52){$abc$}


\put(80,-70){\circle{3.0}} \put(197,-70){\circle{3.0}}
\put(80,10){\circle{3.0}} \put(80,65){\circle{3.0}}
\put(197,-5){\circle{3.0}} \put(198,45){\circle{3.0}}
\put(80,-83){$b_1$}\put(196,-83){$b_2$}

\end{picture}

 The map $t$ has the property that each point $x \in
h^{-1}y$ appears in the multiset $t(y)$ with multiplicity equal to
the size of the parts by which it is labelled, so for the branched
covering of the diagram one has $t(y_1) =
[x_{11},x_{12},x_{13},x_{13},x_{13}],\;$
$t(y_2)=[x_{21},x_{21},x_{22},x_{23},x_{23}]$ and
$t(y_3)=[x_{31},x_{31},x_{31},x_{32},x_{33}].$ We note that, despite
superficial appearances, this example can only be described as an
$n$-branched covering for $n \geq 5.$

\vspace{0.7cm}

 \noindent {\bf \S 4 Frobenius
$n$-homomorphisms and transfer maps}

\medskip
The aim of this section is to introduce the concept of an
$n$-transfer for a ring homomorphism and to study their algebraic
properties as special cases of Frobenius $n$-homomorphisms.

\bigskip
\noindent {\bf Definition 4.1} \hspace{0.3cm} Let $A,B$ be
commutative, associative algebras and $f : A \to B$, a ring
homomorphism, then a linear map $\tau : B \to A$ is an {\bf
$n$-transfer} for $f$ if

\begin{enumerate}
\item[(i)] $\tau$ is a  Frobenius $n$-homomorphism;
\item[(ii)]  $\tau(f(a)b)=a \tau(b)$, that is, $\tau $ is a map of $A$-modules and
\item[(iii)] $f \tau : B \to B$ is the sum of the identity and a Frobenius
$(n-1)$-homomorphism $g :B \to B.$

\end{enumerate}

We denote the linear subspace $\{ b \in B : g(b)=-b\}$ by $L.$

\bigskip
\noindent {\bf Proposition 4.2} \hspace{0.3cm}{\it If $f : A \to B$
is a ring homomorphism, $\tau : B \to A$ is an $n$-transfer for $f$
and $g:B \to B$ is as above then
\begin{enumerate}
\item[(i)] the composition $\tau  f : A \to A $ is multiplication by $n.$
\item[(ii)] there is a split exact sequence of $A$-modules
$$ 0 \to L \to B \stackrel{\tau}{\to} A \to 0.$$
\item[(iii)] $gf = (n-1)f$ and
$\tau g = (n-1) \tau.$

\end{enumerate} }

\bigskip
\noindent {\bf Proof}

Taking $b=1$ in 4.1 (ii) and because $\tau(1)=n$ the result is
immediate.

The map $\tau $ is split by $ f/n.$ If $\tau(b)=0$ then $0 =f
\tau(b) = (1+g)(b)$ by 4.1 (iii) and so the kernel of $\tau$ is
identified with $L.$

The equations  of (iii) are immediate consequences of the
associativity of composition for $f \tau f$ and $\tau f \tau.$

\bigskip
\noindent {\bf Example 4.3} \hspace{0.3cm}
 In the case where $A, \; B$ are
affine algebras, the relations in Proposition 4.2 (iii) can be very
useful. For example, when $A = B = \mathbb{C}[z]$ and $n=2,$ let
$f(z) =p$ and $g(z)=q$ then $gf(z) = p(q(z))$ but this equals $p(z)$
and so $g$ has degree 1.

\bigskip
\noindent {\bf Example 4.4} \hspace{0.3cm} This is another
application of Proposition 4.2 and we only consider the case $n=2.$
The splitting $B = L \oplus A$ is given by identifying the image of
the monomorphism $f :A \to B$ with $A$.  We show that $xy \in A$ for
all $x,y \in L$. From the splitting we obtain $xy=a+z$ where $a \in
A,\;z \in L$; since $g(\ell)=-\ell$ for all $\ell \in L$ and since
$n=2$, 4.2(iii) gives that $g(a)=a$. Using the fact that $g$ is a
ring homomorphism we get that $g(xy)=g(x)g(y)=(-x)(-y)=xy=a+z$ but
also $g(xy)=g(a+z)=a-z$. Hence $z=0,$ showing that $xy \in L.$

\bigskip
\noindent {\bf Proposition 4.5} \hspace{0.3cm} {\it Let $f :A \to B,
\; g:B \to C$ be ring homomorphisms and $\tau:B \to A, \; \sigma :C
\to B$ be $n$- and $m$-transfers for $f, \; g$ respectively. Then
$\sigma \tau : C \to A$ is an $nm$-transfer for the composition
$gf:A \to C.$}

\bigskip \noindent {\bf Proof}

We check the conditions of Definition 4.1. The first condition is
Theorem 2.1 and the other two follow by a direct calculation.

\bigskip
\noindent {\bf Theorem  4.6} \hspace{0.3cm} {\it  If $A,B$ are
algebras with no nilpotent elements and $\tau :B \to A$ is an
$n$-transfer for the ring homomorphism $f :A \to B,$ then
$K_{\tau}=0.$}

\bigskip
{\bf Proof}

Take $b \in K_{\tau}$, then by Gugnin's result, $\tau(b)=\tau(b^2)=
\dots = \tau(b^n)=0.$ By Proposition 4.2 (ii), $g(b^r)=-b^r$ for all
$1 \leq r \leq n.$ But since $g$ is a Frobenius
$(n-1)$-homomorphism, $\Phi_n(g)(b)=0$ and by Definition 1.1 this is
equivalent to the vanishing of the determinant of the matrix

$$
   \left(
\begin{array}{cccccc}
-b & 1 & 0 & 0& \dots &0 \\
-b^2 & -b & 2  & 0& \dots &0 \\
\vdots &\vdots & \vdots& \ddots& & \vdots \\
\vdots&&&&&\vdots\\
-b^{n-1} & -b^{n-2} & -b^{n-3} & \dots & -b
& {n-1 } \\
-b^n &-b^{n-1} & -b^{n-2}  &  \dots &-b^2 &-b

\end{array} \right)
            $$

By adding $b$ times the second column and $b^2$ times the third
column etc to the first column one sees that this determinant equals
$(-1)^n \; n! \; b^n$ and since there  are no nilpotent elements, we
deduce that $b=0.$

\medskip
When $X$ is a compact Hausdorff space, $C(X)$ will denote the algebra
  of continuous functions $X \to \mathbb{C}$ with the supremum norm.

\bigskip
\noindent {\bf Definition 4.7} \hspace{0.3cm} The {\bf direct image}
$t_! : C(X) \to C(Y) $ associated with a continuous map $t : Y \to
{\rm Sym}^n(X)$ is defined by $(t_! \phi)(y) = \sum \phi(x_r)$ and
$t(y) =[x_1,x_2,\ldots, x_n]$.

 \bigskip
\noindent {\bf Theorem 4.8} \hspace{0.3cm} {\it If  $X,Y$ are
compact Hausdorff spaces, then the
  set of all continuous Frobenius $n$-homomorphisms $C(X) \to C(Y)$
  can be identified with the space of continuous maps $Y \to {\rm
  Sym}^n(X).$}

\bigskip
\noindent
{\bf Proof}

The map $t_! $ is the sum of $n$ ring homomorphisms and so is a
Frobenius $n$-homomorphism.

Conversely, suppose that $f :  C(X) \to C(Y)$ be a Frobenius
$n$-homomorphism and  let $\mathcal{E}_y : C(Y) \to \mathbb{C}$ be
evaluation at the point $y \in Y$ then the composition
$\mathcal{E}_y f$ is also  a Frobenius $n$-homomorphism and so, by
Theorem 3.1  of \cite{[BR1]} corresponds to a multi-set $
[x_1,x_2,\ldots, x_n]$ in $X.$ This defines the required map  $t: Y
\to  {\rm Sym}^n(X).$

\bigskip
\noindent {\bf Remark 4.9} \hspace{0.3cm} In the case $A=C(X)$, the
map $\Delta_n$ (see Lemma 1.8) corresponds to the identity map on
${\rm Sym}^n(X)$
 regarded as an $n$-valued map from ${\rm Sym}^n(X)$
to $X.$ If $f : A=C(X) \to B=C(Y)$ is a Frobenius $n$-homomorphism,
the ring homomorphism $\dis{\frac{\Phi_n(f)}{n!}}$ corresponds to $t
: Y \to {\rm Sym}^n(X)$ (see Theorem 1.9).

\bigskip
\noindent {\bf Example 4.10} \hspace{0.3cm}  The linear map  $ C(X)
\to C(X)$ given by $\phi
  \to n \phi$ is a  Frobenius
$n$-homomorphism and corresponds to the diagonal map $X \to {\rm
Sym}^n(X).$

\bigskip

\noindent {\bf Example 4.11} \hspace{0.3cm}  If $s : Y \to {\rm
Sym}^m(X)$ and $t : Y \to {\rm Sym}^n(X)$ give rise to the Frobenius
$m,n$-homomorphisms $s_!,t_! : C(X) \to C(Y)$ then the composition
$$Y \to {\rm Sym}^m(X) \times {\rm Sym}^n(X) \to
{\rm Sym}^{m+n}(X)$$ corresponds to  $s_!+t_! : C(X) \to C(Y).$

\bigskip

\noindent {\bf Example 4.12} \hspace{0.3cm}  If $s:Y \to {\rm
Sym}^n(X)$ and $t : Z \to {\rm Sym}^m(Y)$ are continuous and $s_!
:C(X) \to C(Y), \; t_! : C(Y) \to C(Z)$ are the corresponding $n$-
and $m$-Frobenius homomorphisms, then the composition $$ Z \to {\rm
Sym}^m(Y) \to {\rm Sym}^m{\rm Sym}^n(X) \to{\rm Sym}^{mn}(X) $$
corresponds to the composition $t_!s_!:C(X) \to C(Z)$ which is a
Frobenius $mn$-homomorphism.

\vspace{0.7cm}

 \noindent {\bf \S 5 Frobenius $n$-homomorphisms and $n$-branched coverings }

\medskip

The aim of this section is to characterise $n$-branched coverings in
terms of rings of continuous functions and  Frobenius
$n$-homomorphisms.

\medskip A continuous map $h : X \to Y$ induces a ring homomorphism
$h^{\ast} : C(Y) \to C(X).$ If $h$ is an $n$-branched covering, then
as above we have a direct image map $t_! :C(X) \to C(Y)$ which is a
Frobenius $n$-homomorphism. We consider properties of $t_!$ which
will ensure that $h$ is such a covering.

\medskip

In Example 3.3 of a finite group action on $X$, the third property
of Definition 4.1 becomes very simple : Let $G= \{e=g_1,g_2,\ldots ,
g_n\}$ be the group and $h:X \to Y =X/G$ is the map to the space of
orbits, then the map $h^*t_! :C(X) \to C(X)$ corresponds
geometrically to the map $X \to X \times  \ldots \times X \to X
\times {\rm Sym}^{n-1}(X)$ given by $x \to (x,g_2x,g_3x,\ldots
,g_nx)\to (x,[g_2x,g_3x,\ldots ,g_nx]).$

\medskip
In Theorems 5.1 and 5.2, $X,Y$ will denote compact Hausdorff spaces.

\bigskip
\noindent {\bf Theorem 5.1} \hspace{0.3cm} {\it Given an
$n$-branched covering $h :X \to Y, \;\;t: Y \to {\rm Sym}^n(X)$, the
direct image $t_! : C(X) \to C(Y)$ is an $n$-transfer for the ring
homomorphism $h^* : C(Y) \to C(X).$}

\bigskip\noindent
{\bf Proof}

We check three properties:

 As noted in the proof of Theorem 4.8,
$t_!$ is a Frobenius $n$-homomorphism.

By  definition, one has that for $\phi \in C(X)$, $t_!(\phi)(y)=
\phi(x_1) + \ldots + \phi(x_n)$ where $t(y) = [x_1,x_2,\ldots,
x_n]$. Therefore

$$t_!(h^*(\psi) \phi)(y)= \psi(h(x_1))\phi(x_1) + \ldots +
\psi(h(x_n))\phi(x_n)$$  $$ =\psi(y)\phi(x_1) + \ldots +
\psi(y)\phi(x_n)= \psi(y)t_!(\phi) (y).$$ Hence $t_!(h^*(\psi) \phi)
=\psi(y)t_!(\phi).$

The third property follows immediately from Definition 3.1 (ii).

\medskip

The converse of Theorem 5.1 is

\bigskip
\noindent {\bf Theorem 5.2} \hspace{0.3cm} {\it Given a continuous
map $h : X \to Y$ and  a continuous $n$-transfer $\tau$ for
  $h^{\ast} : C(Y) \to C(X),$ then $h $ is an $n$-branched
  covering.}

\bigskip\noindent
{\bf Proof}

By the above, a continuous $n$-transfer $\tau : C(X) \to C(Y)$
corresponds to a continuous map $ t : Y \to {\rm Sym}^n(X)$ which is
such that $ th : X \to {\rm Sym}^n(X)$ is the diagonal map and $
{\rm Sym}^n(h)t : Y \to {\rm Sym}^n(X) \to {\rm Sym}^n(Y)$ is of the
form $y \to [y_1,y_2,\ldots, y_n]$ with $y_1 =y$.

\bigskip
More generally, we see that the Frobenius $n$-homomorphism $f$
 corresponding to $t: Y \to {\rm Sym}^n(X)$ is the sum of Frobenius
 $n_1$-,$n_2$-,$\ldots,n_k$-homomorphisms $f_1,f_2,\ldots, f_k$ (where
 $n=n_1+n_2+\ldots+n_k$) if and only if $t $ factors as $ Y \to {\rm
 Sym}^{n_1}(X)\times {\rm Sym}^{n_2}(X)\times \ldots \times {\rm
 Sym}^{n_k}(X) \to  {\rm Sym}^n(X)$ where the last map is concatenation.

\bigskip

Using Theorem 3.4 and Corollary 3.6 of \cite{[BR1]} which consider
the relationship between Frobenius $n$-homomorphisms on affine
algebras and symmetric powers of algebraic varieties, one can in a
similar way prove

\bigskip
\noindent {\bf Theorem 5.3 } \hspace{0.3cm} {\it Let $A,B $ be
finitely generated commutative algebras and $f : A \to B$ a ring
homomorphism; let $V(A),V(B)$ be the corresponding varieties and $h:
V(B) \to
 V(A)$ the map corresponding to $f$. Then $h$ is an $n$-fold branched covering if and
only if there is an $n$-transfer $B \to A$, for $h.$}

\bigskip

{\bf Acknowledgements}
\medskip

The research on which this is based was mainly carried out during
visits by VMB to the School of Mathematics, University of Edinburgh
supported by the Engineering and Physical Sciences Research Council.

\pagebreak

\vspace{1cm}

\begin{tabular}{ll}

 Steklov Mathematical Institute, RAS, &
School of Mathematics,\\

 Gubkina 8, &    University of Edinburgh,\\

119991 \hspace{0.1cm} Moscow       & Edinburgh EH9 3JZ\\

{\it and } & {\it and}\\

 School of Mathematics,       &Heilbronn Institute for Mathematical Research, \\
University of Manchester & University of Bristol,\\
Manchester M13 9PL& Bristol BS8 1TW\\
&\\
buchstab@mendeleevo.ru & E.Rees@bristol.ac.uk

\end{tabular}
\end{document}